\documentclass[11pt,a4paper]{amsart}

\usepackage{a4}
\usepackage{amssymb}
\usepackage{latexsym}

\title{Special Lagrangian cones in $\C^3$ and primitive harmonic maps.}
\author{Ian McIntosh}
\address{Department of Mathematics \\ University of York \\ Heslington, York
YO10 5DD, UK}
\email{im7@york.ac.uk}
\subjclass{53C43,58E20}

\newcommand{\SL}{special Lagrangian\ }

\newcommand{\C}{\mathbb{C}}
\newcommand{\R}{\mathbb{R}}
\newcommand{\G}{\Gamma}
\renewcommand{\P}{\mathbb{P}}
\newcommand{\CP}{\mathbb{CP}}
\newcommand{\Z}{\mathbb{Z}}

\newcommand{\End}{\mathrm{End}}
\newcommand{\Hom}{\mathrm{Hom}}

\newcommand{\fa}{\mathfrak{a}}
\newcommand{\fh}{\mathfrak{h}}
\newcommand{\fk}{\mathfrak{k}}

\newcommand{\fg}{\mathfrak{g}}
\newcommand{\fp}{\mathfrak{p}}

\newcommand{\diag}{\mathrm{diag}}
\newcommand{\Ad}{\mathrm{Ad}}

\newcommand{\caA}{\mathcal{A}}
\newcommand{\caE}{\mathcal{E}}
\newcommand{\caL}{\mathcal{L}}
\newcommand{\caO}{\mathcal{O}}

\newtheorem{theor}{Theorem}
\newtheorem{prop}{Proposition}
\newtheorem{lem}{Lemma}

\begin{document}
\begin{abstract}
In this article I show that every special Lagrangian cone in $\C^3$
determines, and is determined by, a primitive harmonic surface in the
6-symmetric space $SU_3/SO_2$. For cones over tori, this allows us to use the
classification theory of harmonic tori to describe the construction of 
all the corresponding
special Lagrangian cones. A parameter count is given for the space of these,
and some of the examples found recently by Joyce are put into this context.
\end{abstract}

\maketitle

\section{Introduction.}

Special Lagrangian submanifolds are the subject of a great deal of attention
at present, largely due to their central role in a number of conjectures
concerning the mirror symmetry of Calabi-Yau 3-folds (see, for example,
\cite{Gro,JoySYZ,StrYZ}). In particular, some effort has been put into the
construction of examples of special Lagrangian cones in $\C^3$
\cite{Has,Joy}, since these
provide local models for the singularities which may occur in special
Lagrangian fibrations of a Calabi-Yau 3-fold. It is not hard to see that a
special Lagrangian cone in $\C^3$ determines, as its link (its intersection
with the unit 5-sphere), a minimal Legendrian
surface in $S^5$ and by Hopf projection a minimal Lagrangian surface in $\CP^2$,
and (locally, at least) the reverse is true. Now, minimal surfaces in spheres
and complex projective spaces possess a well-developed Riemannian twistor
theory \cite{Bur} which asserts that each may be lifted to a harmonic map
into a flag manifold which possesses a sort of horizontal holomorphicity: 
this lift is called a primitive harmonic lift. The first aim of this note is
to describe the relevant twistor theory for special Lagrangian cones in
$\C^3$. 

The second aim of this note is to describe the construction of all minimal
Lagrangian tori in $\CP^2$, which is equivalent to the construction of special Lagrangian
cones over tori. Here I adapt the construction of non-isotropic minimal tori
given in \cite{McI95,McI96} to encode the extra condition
of being Lagrangian. The construction of \cite{McI95,McI96} assigns to each
minimal torus its spectral data $(X,\lambda,\caL)$ consisting of a real 
algebraic curve $X$, a degree three function $\lambda$ on $X$ and a line
bundle $\caL$ over $X$, with $\lambda,\caL$ both respecting the real
involution carried by $X$. The extra
condition of being Lagrangian requires that $X$ possess an extra holomorphic
involution $\mu$ which $\lambda,\caL$ must also respect (this much has also been noted
by \cite{Sha,MaMa}). In particular, 
$\caL$ lives on a translate of the Prym variety $Prym(X,\mu)$.

My interest in this topic grew out of conversations with
Dominic Joyce concerning his construction of special Lagrangian cones over tori
\cite{Joy}. As I became more involved it emerged
that certain aspects of the problem had been tackled, 
from one angle or another, 
at least five times in the last decade. The unifying object 
is the p.d.e
\begin{equation}
\label{eq:Tzitzeica}
\frac{\partial^2u}{\partial z\partial\bar z} = e^{-2u}-e^u,
\end{equation}
which is known either as Tzitz\'{e}ica's equation
or the Bullough-Dodd-Zhiber-Shabat equation. In our context it
occurs as the equation for the metric $e^u|dz|^2$ induced on
a torus by its minimal Lagrangian immersion in $\CP^2$. This is a
soliton equation and its doubly periodic solutions are all of finite gap type:
this fact can be deduced from \cite{BolPW} by recognising (\ref{eq:Tzitzeica})
as the affine $\fa_2^{(2)}$ Toda equations.
It seems to have been first studied in this particular form by Sharipov
\cite{Sha} while he was studying ``zero mean curvature, complex normal''
tori in $S^5\subset \C^3$. These are, one quickly realises, none
other than the minimal Legendrian tori. 
Sharipov studies the finite gap solutions and gives fairly detailed formulae
for the immersion in terms of Riemann $\theta$-functions although he spares no
room for examples. In particular, he points out that the spectral curve must
possess an extra holomorphic involution.
A few years later, Castro \& Urbano \cite{CasU}
solved a reduction of (\ref{eq:Tzitzeica}) to an o.d.e while constructing
$S^1$-equivariant minimal Lagrangian tori in $\CP^2$. Haskins \cite{Has}
performed a similar computation to find $S^1$-equivariant minimal
Legendrian tori in $S^5$. 
Very recently, and following on from Sharipov,
Ma \& Ma \cite{MaMa} gave explicit formulae for the 
minimal Lagrangian maps into $\CP^2$ of finite type 
using Riemann $\theta$-functions, although they do not seem to have noticed 
that their formulae arise from 
Sharipov's by Hopf projection $S^5\to\CP^2$.
Finally, Joyce \cite{Joy} has given equivariant
and non-equivariant minimal Lagrangian tori in $\CP^2$. His immersions are all
given in terms of elliptic functions.

These results can be seen from a broader perspective, for there is
a twistorial picture which will hold for cones over a surface of any
topological type. In section three I describe 
two different but isomorphic primitive harmonic maps which determine a given
special Lagrangian cone 
in $\C^3$.  First, the Gauss map of the
cone itself factors through a harmonic surface which admits a primitive lift
into $SU_3/SO_2$. Second, by applying the twistor theory of \cite{BolPW,Bur}
for the minimal
surfaces in $\CP^2$, we construct a primitive lift into $SU_3/S^1$. 
We view this as an outer 6-symmetric
space, an $S^1$-bundle over the 3-symmetric flag manifold which occurs in
\cite{BolPW,Bur}. The extra property of being Lagrangian imposes
an extra involutive symmetry, corresponding to the outer
involution with which we pass from $\fa_2$ to $\fa_2^{(2)}$. 

In section four we restrict our attention to cones over tori.
It turns out, due to a result of Bolton et al.\
\cite{BolPW}, that every minimal Lagrangian torus in $\CP^2$
lifts to a minimal Legendrian torus in
$S^5$ and thus there is a one-to-one correspondence between minimal Lagrangian
tori in $\CP^2$ and special Lagrangian cones in $\C^3$ whose link is a torus.
I recall how to construct minimal tori in $\CP^2$ and adapt this to
the case of Lagrangian tori: I have tried to do this quickly without leaving out
important details. In short, the spectral curve $X$ of \cite{McI95,McI96} 
must be a double cover $X\to Y$ of another curve, branched only at the unique zero $P_0$
and pole $P_\infty$ of the degree three spectral parameter $\lambda:X\to\hat\C$. This
means $X$ has even genus $2p$, $p=\mathrm{genus}(Y)$, and its Jacobian is isogenous to 
$Jac(Y)\times Prym(X,Y)$ (i.e.\ there is a homomorphism of $Jac(Y)\times Prym(X,Y)$ onto
$Jac(X)$ with finite kernel).
The line bundle $\caL$ lies in a translate of $Prym(X,Y)$ and we can obtain families of
minimal tori by moving $\caL$ along this translate. A parameter count suggests that
minimal Lagrangian tori exist for any spectral genus $2p$: in that case there will be
large (real $(p-2)$-dimensional) 
continuous families of these obtained by moving along the Prym variety.

To finish we consider the case where the spectral curve has genus four; i.e., 
$p=2$. Then $Prym(X,Y)$ is
two dimensional and this forces the solutions of Tzitz\'{e}ica's equation to be doubly
periodic (although this is not sufficient to ensure the minimal map closes into a 2-torus). 
There are two types of spectral curves which can arise: I show 
that for one type $Prym(X,Y)$ is isogenous to a product of elliptic curves. This
means the corresponding minimal maps may be constructed using elliptic functions. 
It is examples of this type which have been found by Joyce \cite{Joy}.

\smallskip\noindent
\textbf{Acknowledgements.} I would like to thank Dominic Joyce for sharing his results 
with me at an early stage and for encouraging me to think about them. I also thank Mark
Haskins for his comments on the first draft.

\smallskip\noindent
\textbf{Notation.} For vectors $v_1,\ldots,v_{n+1}\in\C^{n+1}$ we will use
$v_1\wedge\ldots\wedge v_{n+1}$ to denote the \textit{real} oriented $(n+1)$-plane
with span $Sp_\R\{v_1,\ldots,v_{n+1}\}$. Also, $(v_1,\ldots,v_{n+1})$ will 
denote the matrix with $j$-th column $v_j$. For any Lie group $G$ its
connected component of the identity will be denoted $G^0$. Finally, set 
$\hat\C=\C\cup\{\infty\}$.

\section{Preliminaries.}

A real oriented $(n+1)$-plane $V\subset\C^{n+1}$ is \SL if it lies in the $SU_{n+1}$
orbit of the $(n+1)$-plane $e_1\wedge\ldots\wedge e_{n+1}$ where $\{e_j\}$ is the 
standard unitary basis for $\C^{n+1}$. Let $SL$ denote the set of all these,
then $SL\cong SU_{n+1}/SO_{n+1}$. Since $U_{n+1}\subset SO_{2n+2}$ is the subgroup
of symplectic isometries of $\C^{n+1}$, equipped with its standard K\" ahler
form and Riemannian metric, every element of $SL$ is 
certainly Lagrangian. More generally we say $V$ is \SL with phase $\theta\in[0,2\pi)$ if 
\[
V\in
SL_\theta = \{g(e_1\wedge\ldots\wedge e_{n+1}):g\in U_{n+1},\det(g)=e^{i\theta}\}.
\]

A real oriented $(n+1)$-submanifold $N\subset\C^{n+1}$ is \SL if its Gauss map 
$\gamma_N:N\to Gr_\R(n+1,\C^{n+1})$ takes values in $SL$. The importance of this
property lies in the following two results of Harvey \& Lawson \cite{HarL}.
First, if $N$ is \SL then it is calibrated by the real $(n+1)$-form
\[
\Omega_0 = \mathrm{Re}(dz^1\wedge\ldots\wedge dz^{n+1})
\]
and is therefore both minimal and Lagrangian. Conversely, if $N$ is a 
connected minimal Lagrangian $(n+1)$-fold then it is congruent to a \SL 
$(n+1)$-fold. In particular, every minimal Lagrangian $N$ 
is \SL with phase $\theta$ for some constant $\theta\in[0,2\pi)$.

Let $\hat\pi:\C^{n+1}\setminus\{0\}\to S^{2n+1}$ denote the projection along rays
to the unit sphere $S^{2n+1}\subset\C^{n+1}$. We will say $N\subset\C^{n+1}$ is
a (regular) cone if $N=\hat\pi^{-1}(M)\cup\{0\}$ for some  connected,
oriented submanifold $M\subset S^{2n+1}$ and write $N=C(M)$. It follows 
that $M=C(M)\cap S^{2n+1}$ and $M$ is
called the link of the cone $C(M)$. Since $C(M)$ is usually singular at the origin
we will say $C(M)$ is minimal whenever $C(M)\setminus\{0\}$ is minimal. It is
obvious that $C(M)$ is Lagrangian whenever $M$ is Legendrian for the contact 
structure $S^{2n+1}$ inherits from $\C^{n+1}$. Moreover, since $\hat\pi$ is
a Riemannian submersion with geodesic fibres, $C(M)$ is minimal if and only if
$M$ is minimal. We deduce (as did Haskins \cite{Has}) that there is a one-to-one
correspondence between minimal Lagrangian $(n+1)$-cones in $\C^{n+1}$ and minimal
Legendrian $n$-folds in $S^{2n+1}$.

Now let $\pi:S^{2n+1}\to\CP^n$ be the Hopf fibration. Again, this is a Riemannian
submersion with geodesic fibres, so if $M\subset S^{2n+1}$ is minimal then 
so is $\pi(M)$. Further, the pullback $\pi^*\omega$ of the Fubini-Study K\" ahler form 
is, up to a constant scalar, the curvature for the natural metric connection carried by 
the Hopf bundle. Moreover, the horizontal distribution 
is the standard contact structure on $S^{2n+1}$; i.e., 
a Legendrian submanifold of $S^{2n+1}$ is horizontal for $\pi$. 
Therefore every minimal Legendrian immersion $f:M\to S^{2n+1}$ induces a minimal
Lagrangian immersion $\varphi=\pi\circ f:M\to\CP^n$. Conversely, a minimal
Lagrangian $\varphi:M\to\CP^n$ has $\varphi^*\omega = 0$ and therefore
has a horizontal 
lift of its universal cover $f:\tilde M\to S^{2n+1}$. It is a straightforward
calculation to show that: (i) since $f$ is horizontal the tension fields satisfy
$\pi_*\tau(f)=\tau(\varphi)$ (see, for example, \cite[p 20]{EelL}), and; (ii) 
$\tau(f)$ is horizontal (using \cite[\S 2]{ONe}). 
Therefore $f$ is minimal and, of course, Legendrian. So
for simply connected manifolds there is, up to congruence, 
a one-to-one correspondence between
minimal Lagrangian immersions in $\CP^n$ and minimal Legendrian immersions 
in $S^{2n+1}$.

\section{Twistor theory for \SL cones in $\C^3$.}

Let us now fix a \SL cone $C(M)\subset\C^3$ with link $f:M\to S^5$ and the
corresponding minimal Lagrangian $\varphi:M\to\CP^2$. Our aim is to show that
any one of these maps is characterised by the existence of a primitive 
harmonic lift of $M$ into a certain 6-symmetric space. Let us first recall
(from \cite[\S 1]{Bur}) the idea of a primitive harmonic map into a $k$-symmetric 
space.

A $k$-symmetric space is a reductive homogeneous space $G/K$ for a compact 
semisimple Lie group $G$ in which the isotropy group $K$ is the subgroup
of fixed points of an order $k$ automorphism $\sigma$ on $G$ (or, more generally,
we can take $K$ to be a closed subgroup of the fixed point group which shares 
its Lie algebra). This
automorphism induces an automorphism on $\fg$, the Lie algebra of $G$,
which we will also call $\sigma$, with
subalgebra of fixed points $\fk$, the Lie algebra of $K$. The complexification
$\fg^\C$ decomposes into the direct sum of
eigenspaces $\fg_j$ for eigenvalues $\epsilon^j$ where $\epsilon$ is a
primitive $k$-th root of unity. The reductive decomposition
of $\fg$ is characterised by
\[
\fg=\fk +\fp;\quad \fk^\C = \fg_0,\ \fp^\C=\sum_{j\neq 0}\fg_j.
\]  
Clearly $K$ acts by adjoint action on each $\fg_j$. 
We recall that
$T(G/K)$ is naturally isomorphic to $[\fp]=G\times_K\fp$
(the quotient of $G\times\fp$ by right $K$-action). We can identify this with
a subbundle of $G/K\times \fg$ so that fibrewise the isomorphism is 
\[
[\fp]_q\to T_q(G/K);\ X\mapsto \frac{d}{dt}(e^{tX}\cdot q)|_{t=0}.
\]
In particular $[\fp]_q=\Ad g\cdot \fp$ whenever $q=gK$.
The complexification
$T^\C(G/K)$ possesses subbundles $[\fg_j]$ for $j\neq 0$. For $k\geq 3$, a smoothly immersed 
surface $\psi:M\to G/K$ is said to be primitive whenever 
\begin{equation}
\label{fhol}
d\psi^{1,0}:TM^{1,0}\to [\fg_1].
\end{equation}
This is a weakened version of the $J$-holomorphic condition which one finds in Riemannian
twistor theory. Let $F\in\End(T^\C G/K)$ be the bundle endomorphism satisfying $F^3=-F$ 
and having $[\fg_1]$ as its bundle of $i$-eigenspaces, $[\fg_{k-1}]$ as its bundle
of $-i$-eigenspaces and kernel $\oplus_{j=2}^{k-2}[\fg_j]$. Then $F$ is
one of the horizontal $f$-structures studied by Black \cite[\S 6]{Bla}.
The condition \eqref{fhol} means $\psi$ is $f$-holomorphic with respect to this $f$-structure.
By results of Black \cite[Thm 4.3, Lemma 5.3]{Bla} we have:
(1) every primitive map is harmonic (for any $G$-invariant metric on $G/K$ for which $[\fg_1]$
is isotropic);
(2) if $H\subset G$ is any closed subgroup containing $K$ then post-composition of
any primitive map $\psi$ with the homogeneous projection $G/K\to G/H$ again yields
a harmonic map. 

Given the \SL cone $C(M)$ we will find two primitive harmonic maps. The first,
$\psi_1:M\to SU_3/SO_2$, arises from the Gauss map of $C(M)$. The second,
$\psi_2:M\to SU_3/S^1$, arises by considering how the Lagrangian condition
interacts with the primitive lift assigned by
\cite{BolPW,Bur} to any minimal surface in $\CP^2$. We will show that 
these are basically the same map. 

\subsection{The Gauss map of $C(M)$.}

Let $\gamma_{C(M)}:C(M)\setminus\{0\}\to SL$ be the Gauss map of $C(M)$. Note that since $C(M)$
is minimal  $\gamma_{C(M)}$ is harmonic by the Ruh-Vilms theorem (see, for example,
\cite[p 18]{EelL}). Since $C(M)$ is a cone $\gamma_{C(M)}$ factors through $M$. Indeed, it
is easy to see that if $\gamma_M:M\to Gr_\R(2,\C^3)$ denotes the Gauss map
of the link $f:M\to S^5$ then $\gamma_{C(M)}$ is projection onto $M$
followed by
\begin{equation}
\label{eq:gauss}
f\wedge \gamma_M:M\to SL\cong SU_3/SO_3.
\end{equation}
This too is harmonic since the projection is Riemannian with geodesic fibres.

Now set
\[
FL_1 = \{(v,V)\in S^5\times SL: v\in V\}= \{(ge_1,g(e_1\wedge e_2\wedge e_3)):g\in SU_3\}.
\]
Clearly $FL_1\cong SU_3/K_1$ where $K_1=\{g\in SO_3: ge_1=e_1\}\cong SO_2$ and 
the two natural projections
\[
\pi_1:SU_3/K_1\to SU_3/SU_2\cong S^5;\ \pi_2:SU_3/K_1\to SU_3/SO_3\cong SL
\]
are homogeneous projections.
The map (\ref{eq:gauss}) obviously possesses a lift 
\[
\psi_1:M\to FL_1;\quad \psi_1 = (f,f\wedge\gamma_M).
\]
We can view $FL_1$ as a $6$-symmetric space as follows. Let 
$\hat\sigma\in Aut(SU_3)$ be given by
\[
\hat\sigma(g)=Rg^{-1t}R^{-1},
\]
where $R\in SO_3$ is the rotation about the $e_1$ axis through
angle $\pi/3$. Then $\hat\sigma$
is an outer automorphism of order six  with subgroup of fixed
points $K_1$. 
\begin{prop}
\label{prop:psi1}
The map $\psi_1:M\to FL_1$ is primitive harmonic. 
Conversely, any such map projects by
$\pi_1$ to the Gauss map of an immersed \SL cone $\kappa:C(M)\to\C^3$ 
and by $\pi_2$
to a minimal Legendrian map $f:M\to S^5$.
\end{prop}
The proof of this will follow once we have proved a similar result for
the primitive lift of the minimal Lagrangian map $\varphi:M\to\CP^2$.

\subsection{Minimal Lagrangian surfaces in $\CP^2$.}

As we noted above, every minimal Lagrangian $\varphi:M\to\CP^2$ lifts
(locally or on the universal cover) horizontally to a minimal Legendrian
map into $S^5$. For simplicity let us assume this lift is global $f:M\to
S^5$; i.e., $\varphi$ arises from a \SL cone in $\C^3$. We will use ideas
from \cite{BolW,Bur} to lift $\varphi$ into the manifold 
\[
FL_2 = \{(w,W)\in S^5\times Gr_\C(2,\C^3):w\in W\}.
\]
The motivation for this is that \cite{BolPW,Bur} have already shown that
any minimal surface in $\CP^2$ lifts to a primitive map into the flag manifold
\[
Fl_2 = \{(W_1,W_2)\in \CP^2\times Gr_\C(2,\C^3):W_1\subset W_2\}.
\]
We will show that the extra symmetry of being Lagrangian allows us to lift
to the $S^1$-bundle $pr_2:FL_2\to Fl_2$ and the lift will also be primitive.
If we think of $Fl_2$ as the $3$-symmetric space for the Coxeter-Killing
automorphism $\nu$ (whose definition will be recalled shortly) then $FL_2$ is the
$6$-symmetric space for an automorphism $\sigma$ which is the product of $\nu$ and
a commuting involution $\mu$.

To begin, let $L\to\CP^2$ denote the tauotogical line bundle, then $\varphi$ determines (and
is determined by) $\ell_0=\varphi^*L \subset M\times\C^3$. Any smooth subbundle $\ell\subset M\times\C^3$
can be equipped
with a holomorphic structure for which a local section is holomorphic
whenever $\pi_\ell(\partial s/\partial\bar z)=0$, where $\pi_\ell$ is the 
orthogonal projection onto $\ell$
and $z$ is a local complex coordinate
on $M$.  Recall that $\varphi$ is harmonic precisely when the map
\[
d\varphi(\frac{\partial}{\partial z}):\ell_0\to\ell_0^\perp;\ s\mapsto 
\pi_{\ell_0}^\perp(\frac{\partial s}{\partial z}),
\]
is holomorphic. It is well known that in this case
the new line bundles $\ell_1,\ell_2,\ldots$ defined by 
\begin{equation}
\label{eq:sequence}
\ell_{j+1} = \pi_{\ell_j}^\perp(\frac{\partial}{\partial z}\ell_j),
\end{equation}
each determine a harmonic map $\varphi_j:M\to\CP^n$ for which
$\ell_j=\varphi_j^*L$. The sequence
$\ell_0,\ell_1,\ldots$ is called the harmonic sequence (see eg.\ \cite{BolW}).  
The extra condition of conformality is equivalent to $\ell_2\perp\ell_0$.
Thus from a minimal surface in $\CP^2$ 
we obtain an orthogonal harmonic sequence $\ell_0,\,\ell_1,\ell_2$ which
decomposes $M\times \C^3$. The primitive harmonic lift of $\varphi$ is the map
\[
\psi:M\to Fl_2;\ \psi = (\ell_0,\ell_0\oplus\ell_1).
\]
Now, since $\varphi$ is minimal Lagrangian with the lift $f:M\to S^5$, we can define
\[
\psi_2:M\to FL_2;\ \psi_2 = (f,\ell_0\oplus\ell_1).
\]
We will show that this is primitive harmonic with respect to a particular $6$-symmetric
space structure on $FL_2$. 

First, notice that $FL_2$ is the $SU_3$-orbit of $(w_0,W_0)$, where
\[
w_0=e_1,\ W_0= Sp_\C\{e_1,e_2\}.
\]
The isotropy group is $K_2=\{\diag(1,a,a^{-1})\in SU_3\}$. Therefore 
$FL_2\cong SU_3/K_2$ and it fibres over the flag manifold $Fl_2\cong
SU_3/T$ (where $T$ is the maximal torus of diagonal matrices) with 
homogeneous projection
\[
pr_2:FL_2\to Fl_2;\ (w,W)\mapsto (Sp_\C\{w\},W).
\]
We view $FL_2$ as the $6$-symmetric space corresponding
to the order six automorphism $\sigma\in Aut(SU_3)$ given by the product
$\sigma=\mu\nu$ where:
(i) $\nu$ is the Coxeter-Killing automorphism $\nu(g)=SgS^{-1}$ where
$S=\diag(1,\epsilon,\epsilon^2)$ for $\epsilon=\exp(2\pi i/3)$; (ii)
$\mu$ is the involution $\mu(g) = Tg^{-1t}T^{-1}$ where 
\begin{equation}
\label{eq:T}
Te_1=e_1,\ Te_2=e_3,\ Te_3=e_2.
\end{equation}
These two automorphisms commute and so $\sigma$ is well-defined, with subgroup
of fixed points $K_2$. On $\fg=\mathfrak{su}_3$ we take the eigenvalues of $\sigma$
to be $\{(-\epsilon)^j:j=0,\ldots 5\}$: this makes it 
easier to work with the eigenspaces $\fg_j$ of $\sigma$, by observing that
\begin{equation}
\label{eq:grading}
\fg_j = \fg^\mu_{j\bmod 2}\cap \fg^\nu_{j\bmod 3}.
\end{equation}
In particular
\[
\fg_1=\{X\in\fg^\C: X=
\begin{pmatrix}
0 & 0 & \alpha\\
\alpha & 0 & 0\\
0 & \beta & 0
\end{pmatrix},\alpha,\beta\in\C\}
\]
Let $\fg=\fk_2\oplus \fp_2$ be the corresponding reductive splitting and make
the usual identification
\begin{equation}
\label{eq:tangent}
T^\C FL_2\cong [\fp_2^\C]=SU_3\times_{K_2} \fp_2^\C.
\end{equation}
The primitive distribution is the subbundle
$[\fg_1]$ whose fibres can be described as
\[
[\fg_1]_q = \Ad g\cdot\fg_1\subset\fg^\C;\ \hbox{at}\ q=gK_2.
\]
\begin{prop}
\label{prop:psi2}
Let $\varphi:M\to\CP^2$ be minimal Lagrangian with horizontal lift $f:M\to S^5$.
Then we can choose $f$ so that the
lift $\psi_2:M\to FL_2$ is primitive harmonic; i.e., 
$d\psi_2^{1,0}\in\Omega_M^{1,0}\otimes \psi_2^*[\fg_1]$ for the
$6$-symmetric space structure described above.
Conversely, any primitive harmonic $\psi_2$ projects
by $pr_1:FL_2\to S^5$ to a minimal Legendrian surface in $S^5$
and thence by the Hopf map to a minimal Lagrangian surface in $\CP^2$.
\end{prop}
\noindent
\textit{Remark.} This  special choice of horizontal lift $f$ is an artefact of our definition
of $FL_2$. In the next subsection I will show that any choice of $f$ gives rise to some
primitive lift into an isomorphic $6$-symmetric space.
\begin{proof}
We can move $f$ by an isometry of the form $f\mapsto e^{i\theta}f$ 
so as to ensure that the Lagrangian 
map $f\wedge\gamma_M$ is special with phase $3\pi/2$. In terms of a local complex
conformal coordinate $z=x+iy$ we have  
\[
\det(f,\frac{1}{|f_x|}f_x,\frac{1}{|f_y|}f_y)=-i.
\]
Now we can choose a local frame $F$ for $\psi_2$ by taking
\begin{equation}
\label{eq:frame}
F=(f,\frac{1}{|f_z|}f_z,\frac{-1}{|f_{\bar z}|}f_{\bar z}) = (f_0,f_1,f_2).
\end{equation}
Now, $f$ is a horizontal (and therefore holomorphic) section
of $\ell_0$, and since $\varphi$ is conformal we have
\[
\langle f_z,f\rangle =0=\langle f_{\bar z},f\rangle,\ \langle f_z,f_{\bar z}\rangle =
0.
\]
Therefore each $f_j$ is a local section of $\ell_j$ and $F$ is unitary. Further,
\begin{equation}
\label{eq:U}
(f,\frac{1}{|f_x|}f_x,\frac{1}{|f_y|}f_y)=FU\ \hbox{for}\ U=
\begin{pmatrix} 1 & 0 & 0 \\
0 & 1/\sqrt{2} & -i/\sqrt{2}\\
0 & -1/\sqrt{2}& -i/\sqrt{2}
\end{pmatrix}.
\end{equation}
Since $\det(U)=-i$ it follows that $F$ is special unitary. Now recall that using the
identification (\ref{eq:tangent}) we identify
\[
d\psi_2(\partial/ \partial z)= \Ad F\cdot 
(F^{-1}\frac{\partial F}{\partial z})_{\fp_2^\C}.
\]
A straightforward calculation gives
\begin{equation}
\label{eq:dpsi2}
(F^{-1}\frac{\partial F}{\partial z})_{\fp_2^\C}=
\begin{pmatrix}
0 & 0 & |f_{\bar z}|\\
|f_z| & 0 & 0 \\
0 & Q|f_z|^{-2}& 0
\end{pmatrix}
\end{equation}
where $Q=\langle f_{zzz},f\rangle$.
But now we observe that $\langle f_{\bar z z},f\rangle = \langle f_{z\bar z},f\rangle$
implies that $|f_z|=|f_{\bar z}|$ and therefore $d\psi_2^{1,0}$ takes values in the
primitive distribution.
The converse statement in the proposition
follows at once from the result due to Black \cite{Bla} mentioned earlier.
\end{proof}
A standard argument (see, for example, \cite[p 125]{BolPW}) shows
that the cubic form $\det(d\psi_2^{1,0})=Qdz^3$ is
holomorphic on $M$, since $\psi_2$ is harmonic for the $SU_3$-invariant metric. 
Indeed, $Qdz^3$ is globally defined
even if the lift $f:M\to S^5$ is not
global since it is easy to see that $d\psi_2$ is the pullback of $d\psi$ under
$FL_2\to Fl_2$, and $\psi:M\to Fl_2$ exists as a global lift for any minimal surface in
$\CP^2$. The form $Qdz^3$ is the Hopf differential of $\varphi$: its zeroes represent
the higher order singularities of $\psi$.

\subsection{The isomorphism between $FL_1$ and $FL_2$.}

In the previous section we found a lift $\psi_2$ of the minimal Lagrangian map 
$\varphi$ which corresponded to a
\SL cone with phase $3\pi/2$. This phase factor can be avoided if we
perform the following sleight of hand. We identify $FL_2$ with $N/K_2$ where
\[
N=\{g\in U_3:\det(g)=i\}=\{U^{-1}g:g\in SU_3\},
\]
and $U$ is the matrix defined in (\ref{eq:U}). The implicit identification
$N/K_2\cong SU_3/K_2$ allows us to identify $T^\C N/K_2$ with $[\fp_2^\C]$ so that the
primitive distribution has fibre
\[
[\fg_1]_q = \Ad g\cdot \fg_1\subset\fg^\C,\ \hbox{at}\ q=U^{-1}gK_2.
\]
Now if we take $f:M\to S^5$ to be a special Legendrian lift of $\varphi$, with phase
zero, the local frame $F$ in (\ref{eq:frame}) has $\det(F)=i$ so that
$\psi_2=(f,\ell_0\oplus\ell_1)$ is a primitive harmonic lift. 

With this slightly altered point of view
there is a natural isomorphism $\Phi:FL_2\cong FL_1$ for which
$\Phi_*[\fg_1]$ is the primitive distribution for $FL_1$ and for which $\psi_1=
\Phi\circ\psi_2$. Define 
\[
\Phi:FL_2\to FL_1;\ (w,W)\mapsto (v,V),
\]
by $v=w$ and $V=v_0\wedge v_1\wedge v_2$ where
\[
(v,v_1,v_2)=(w,\frac{1}{\sqrt{2}}(w_1-w_2),\frac{-i}{\sqrt{2}}(w_1+w_2))=(w,w_1,w_2)U.
\]
Here $w,w_1,w_2$ is any unitary basis
of $\C^3$ adapted to the flag $Sp_\C\{w\}\subset W\subset\C^3$ 
and with $\det(w,w_1,w_2)=i$. It is easy to check that
$\Phi$ is independent of which basis we choose and 
\[
\det(v,v_1,v_2)=1.
\]
Further, since $U$ is unitary the real $3$-plane $V$ is Lagrangian,
whence $V\in SL$. 
\begin{prop}
The isomorphism $\Phi:FL_2\to FL_1$ has the properties: (i)
$\Phi_*[\fg_1]=[\hat\fg_1]$ where $[\hat\fg_1]$ is the primitive distribution for 
the $6$-symmetric space structure on $FL_1$ determined by $\hat\sigma$ given earlier; 
(ii) $\psi_1=\Phi\circ\psi_2$. Therefore
$\psi_1$ is primitive harmonic whenever $\psi_2$ is.
\end{prop}
\begin{proof}
Clearly $\psi_1=\Phi\circ\psi_2$ by (\ref{eq:U}).
Now let us consider $\Phi_*:T^\C FL_2\to T^\C FL_1$. 
The identification of $FL_1$ with $SU_3/K_1$ gives the usual isomorphism
\[
T^\C_p FL_1\cong \Ad g\cdot \fp^\C_1\subset\fg^\C,\ \hbox{at}\ p=gK_1.
\]
We can represent $\Phi$ by 
\[
\Phi:N/K_2\to SU_3/K_1;\ U^{-1}gK_2\to U^{-1}gUK_2,\ \det(g)=1,
\]
(the reader can check that for any $k\in K_2$ we have $U^{-1}kU\in K_1$). 
Therefore
\[
\Phi_*:[\fp_2^\C]\to [\fp_1^\C];\ \Phi_*(X) = \Ad U^{-1}\cdot X.
\]
Now consider $X\in [\fg_1]_q$. We can write $X=\Ad g\cdot Y$ where $Y\in\fg_1$; i.e.,
\[
-\Ad TS\cdot Y^t= -\epsilon Y.
\]
Then $\Phi_*(X) = \Ad U^{-1}g\cdot Y$ so 
\[
\Phi_*(X)=\Ad U^{-1}gU\cdot \hat Y,\ \hbox{for}\ \hat Y=\Ad U^{-1}\cdot Y.
\]
Therefore
\[
\Ad TS\cdot(\Ad U\cdot \hat Y)^t = \epsilon\Ad U\cdot\hat Y
\]
hence $-\Ad R\cdot\hat Y^t = -\epsilon\hat Y$, where
\[
R = U^{-1}TSU^{-1t} = 
\begin{pmatrix}
1& 0 & 0\\
0 & \cos(\pi/3) & \sin(\pi/3)\\
0 & -\sin(\pi/3) & \cos(\pi/3)
\end{pmatrix}.
\]
But $FL_1$ has the $6$-symmetric space structure corresponding to
the automorphism $\hat\sigma(g)=Rg^{-1t}R$ of $SU_3$, for which the primitive
distribution is determined by
\[
\hat\fg_1 = \{X\in\fg^\C:-\Ad R\cdot X^t = -\epsilon X\}.
\]
Thus $\Phi_*[\fg_1]=[\hat\fg_1]$.
\end{proof}
Proposition \ref{prop:psi1} follows now from proposition \ref{prop:psi2}.

\section{Minimal Lagrangian tori in $\CP^2$.}

Let us now suppose we have a minimal Lagrangian torus $\varphi:T^2\to\CP^2$. 
In this case the Hopf differential $Q$ must be constant. Either $Q=0$, in which case
$\varphi$ is superminimal, or $Q\neq 0$ and $\varphi$ is
superconformal. One knows from \cite[Thm 3.6]{BolW} that a superminimal Lagrangian
surface must have image congruent to $\R\P^2$: the corresponding \SL cone is
$\R^3\subset\C^3$. So we may assume $\varphi$ is superconformal
and therefore possesses a \textit{Toda frame}. This means (see \cite[\S 2]{BolPW})
there is a frame $F:\R^2\to SU_3$ satisfying, in an appropriate 
coordinate $z$, the equation
\begin{equation}
\label{eq:todaframe}
F^{-1}\frac{\partial F}{\partial z} = 
\begin{pmatrix} s_0^{-1}\partial s_0/\partial z & 0 & s_0s_2^{-1}\\
s_1s_0^{-1} & s_1^{-1}\partial s_1/\partial z & 0\\
0 & s_2s_1^{-1} & s_2^{-1}\partial s_2/\partial z
\end{pmatrix},
\end{equation}
where $s_j:T^2\to \R^+$ satisfy the affine $\fa_2$ Toda equations
\begin{equation}
\label{eq:toda}
\frac{\partial^2}{\partial z\partial\bar z}\log (s_j^2) = s_{j+1}^2s_j^{-2} -
s_j^2s_{j-1}^{-2},\ j\in\Z_3.
\end{equation}
Recall that the latter are simply the integrability conditions for the existence
of $F$.
Moreover, it was shown in \cite[Thm 2.5]{BolPW} that Toda frames are essentially 
uniquely determined by this property and
they are also double periodic, to wit, $F:\tilde T^2\to
SU_3$ where $\tilde T^2\to T^2$ is a covering torus with covering group either $1$ or
$\Z_3$. Using these facts let us now prove:
\begin{prop}
Every minimal Lagrangian torus $\varphi:T^2\to\CP^2$ has a doubly periodic primitive
lift $\psi_2:\tilde T^2\to FL_2$ where $\tilde T^2\to T^2$ has covering group either $1$
or $\Z_3$. Consequently $\varphi$ determines a minimal Legendrian $f:\tilde T^2\to S^5$
and an immersed special Lagrangian cone $\kappa:C(\tilde T^2)\to\C^3$.
\end{prop}
\begin{proof}
It suffices to show that, with an appropriate choice of complex coordinate $z$ on $\R^2$
the frame $F=(f_0,f_1,f_2)$ given by (\ref{eq:frame})
satisfies the Toda frame equations. Since the Hopf differential $Qdz^3$ is globally
holomorphic on $T^2$ we can choose $z$ such that $Q=\langle f_{zzz},f\rangle = 1$. Now
define $s_j:\R^2\to\R^+$ by $s_0=1,s_1=|f_z|,s_2=1/|f_z|$. 
The equation (\ref{eq:dpsi2})
shows that $F$ satisfies (\ref{eq:todaframe}) provided we can show that
\[
\langle \frac{\partial f_j}{\partial z},f_j\rangle = 
s_j^{-1}\frac{\partial}{\partial z}s_j,\ j\in\Z_3.
\]
The reader is left to check that these equations are true using the identities
$f_1=s_1^{-1}f_z,f_2 = -s_2f_{\bar z}$ and 
\[
\langle f_z,f\rangle =0,\ \langle f_{z\bar z},f_z\rangle = 0 = \langle f_{z\bar z},f_{\bar
z}\rangle.
\]
The last pair of equations follows from the fact that 
$f_z= d\varphi(f\otimes\partial/\partial
z)$ is a holomorphic section of
$\ell_1$ while $f_{\bar z}= d\varphi(f\otimes\partial/\partial \bar z)$ is an 
anti-holomorphic section of $\ell_2$, since $\varphi$ is harmonic. Here we are
interpreting $d\varphi$ as a 1-form with values in 
$\varphi^{-1}T\CP^2\cong \Hom(\ell_0,\ell_0^\perp)$
(cf.\ \cite{BolW}).
\end{proof}

\subsection{The spectral data for a minimal Lagrangian torus.}

Recall from \cite{BolPW,Bur} that for every minimal torus $\varphi:T^2\to\CP^2$ its primitive lift
$\psi:T^2\to Fl_2$ is of semisimple finite type. This ultimately means that $\varphi$ is
determined, uniquely up to congruence, by its spectral data $(X,\lambda,\caL)$ \cite{McI95}.
Here I will describe what this data is for a minimal Lagrangian torus: basically, the extra
condition of being Lagrangian equips $X$ with a holomorphic involution which $\lambda$ and
$\caL$ must also respect.

As above, let $F=(f_0,f_1,f_2)$ be the Toda frame for a minimal Lagrangian
torus $\varphi:M\to\CP^2$. Set $\alpha=F^{-1}dF$, and write
$\alpha=\alpha_\fk+\alpha_\fp$ for the components in the reductive
decomposition of $\fg$. Since $s_0=1,s_2=s_1^{-1}$ these components satisfy
\[
\alpha_\fk(\frac{\partial}{\partial z}) =
\begin{pmatrix} 
0 & 0 & 0\\ 
0 & s_1^{-1}\partial s_1/\partial z & 0 \\
0 & 0 & -s_1^{-1}\partial s_1/\partial z
\end{pmatrix},\ 
\alpha_\fp(\frac{\partial}{\partial z})=
\begin{pmatrix}
0 & 0 & s_1\\
s_1 & 0 & 0\\
0 & s_1^{-2}& 0
\end{pmatrix}.
\]
Further, $s_1:T^2\to\R^+$ satisfies the affine $\fa_2^{(2)}$ Toda equations
\begin{equation}
\label{eq:tzitz}
\frac{\partial^2}{\partial z\partial\bar z}\log(s_1^2) = s_1^{-4}-s_1^2.
\end{equation}
Notice that $s_1= \parallel \partial\varphi/\partial z\parallel$ 
and therefore $s_1^2$
represents the conformal factor in the metric induced on $T^2$ by 
the immersion $\varphi$ (cf.\ \cite[\S 2]{CasU}).

Now we recall from \cite[\S 3]{BurP} that the lift $\psi_2$ is harmonic
if and only if any frame can be extended to an $S^1$-family (or loop) of frames
\[
F_\zeta:\R^2\to\Lambda^\sigma SU_3 = \{g_\zeta:S^1\to SU_3|\sigma(g_\zeta)=g_{-\epsilon\zeta}\}.
\]
This loop of frames is determined, up to initial condition, by
\[
F_\zeta^{-1}dF_\zeta = \alpha_\zeta=\zeta\alpha_\fp^{1,0} +\alpha_\fk+\zeta^{-1}\alpha_\fp^{0,1}.
\]
We note that $\alpha_\zeta$ actually possesses three symmetries
\[
\nu(\alpha_\zeta)=\alpha_{\epsilon\zeta},\quad  \mu(\alpha_\zeta)=\alpha_{-\zeta},\quad
-\alpha_\zeta^\dagger = \alpha_{\bar\zeta^{-1}},
\]
(the first two can be deduced from (\ref{eq:grading})), which obliges $F_\zeta$ to have the
same symmetries (for every initial condition with these symmetries).

We recall from \cite[p 240]{McI98} that there is a commutative algebra $\caA$ of
complexified polynomial Killing fields
\[
\xi_\zeta:\R^2\to\Lambda^\nu\mathfrak{gl}_3=
\{X_\zeta:S^1\to\mathfrak{gl}_3(\C)|\nu(X_\zeta)=X_{\epsilon\zeta}\},\quad 
d\xi_\zeta=[\xi_\zeta,\alpha_\zeta],
\]
where $\xi_\zeta$ is a Laurent polynomial in $\zeta$ (of bounded degree over $\R^2$). 
By the symmetries of $\alpha_\zeta$, this algebra admits two automorphisms
\[
\mu_*(\xi_\zeta) = \mu(\xi_{-\zeta}),\ \rho_*(\xi_\zeta) = -\bar\xi_{\bar\zeta^{-1}}^t.
\]
The algebra $\caA$ determines the spectral curve $X$ of $\varphi$, as the smooth completion of the
the affine curve $Spec(\caA)$, together with a degree three cover 
$\lambda:X\to\hat\C$
with a single zero $P_0$ and pole $P_\infty$. This cover is the dual of $\C[\zeta^3I]\subset\caA$,
where $I$ is the identity matrix. The real involution $\rho_*$ equips $X$ with a real
involution $\rho$ covering $\lambda\mapsto\bar\lambda^{-1}$, for which $\rho$ fixes 
every point over $|\lambda|=1$, while $\mu_*$ induces a holomorphic
involution $\mu$ on $X$ which fixes only the points $P_0,P_\infty$ (since it must cover
$\lambda\mapsto -\lambda$). Inside $\caA$ we have the subalgebra
$\caA^\mu=\{\xi_\zeta:\mu_*(\xi_\zeta)=\xi_\zeta\}$, which we might naturally think of as the
polynomial Killing fields for $\psi_2$. 
\begin{lem}
$\caA = \caA^\mu[\zeta^3I]\cong \caA^\mu[t]/(t^2-\zeta^6I)$, where $(t^2-\zeta^6I)$ denotes the
ideal generated by this quadratic. 
Therefore the inclusion $\caA^\mu\subset\caA$ is dual to a double cover $\delta:X\to Y$, 
ramified at $P_0$ and $P_\infty$ only, where
$Y$ is the smooth completion of $Spec(\caA^\mu)$. The involution $\mu$ is the swapping of sheets
for this cover. It follows that $X$ has even genus $g=2p$ where $p$ is the genus of $Y$.
\end{lem}
\begin{proof}
We can write each $\xi_\zeta$ as $A_\zeta +\zeta^3 B_\zeta$ where
\[
A_\zeta = \xi_\zeta +\mu_*(\xi_\zeta),\ B_\zeta = \zeta^{-3}(\xi_\zeta -\mu_*(\xi_\zeta)).
\]
Then $A_\zeta,B_\zeta\in\caA^\mu$ and we notice that $\C[\zeta^6I]\subset\caA^\mu$.
It follows that $\caA = \caA^\mu[\zeta^3I]$. Therefore $X$ is the Riemann 
surface for the quadratic extension $t^2-\zeta^6I=0$ over the smooth
completion $Y$ of $Spec(\caA^\mu)$. This is clearly only ramified at
$P_0,P_\infty$. 
\end{proof}

Next we recall from \cite{McI98} that any minimal torus $\varphi$ is determined, 
up to congruence, by the data $(X,\lambda,\caL)$ where $\caL$ is a holomorphic 
line bundle over $X$ of degree $g+2$. To be precise, we identify
$T^2$ with $\C/\Lambda$ for some lattice $\Lambda$ and assign to $\varphi(z)$ the data
$(X,\lambda,\caL_z)$: this whole family is determined by the initial data
$(X,\lambda,\caL_0)$. Let us recall how this works. We will think of
$\caL_z$ as the
dual to the eigenline bundle $\caE_z$ over $X$, defined as follows. First, to each
$\xi_\zeta\in\caA$ we assign an untwisted polynomial Killing field
\[
\hat\xi_{\zeta^3}=\Ad\kappa_\zeta\cdot\xi_\zeta,\quad \kappa_\zeta=\diag(1,\zeta^{-1},\zeta^{-2}).
\]
Since $\Ad\kappa_{\epsilon\zeta}=\Ad\kappa_\zeta\circ\nu$ we deduce that this untwisted loop
really is a function of $\lambda=\zeta^3$. This provides an isomorphic algebra
$\hat\caA=\{\hat\xi_\lambda:\xi_\zeta\in\caA\}$. Notice that this untwisting has not 
changed the symmetry $\mu_*:\hat\caA\to\hat\caA$; 
$\mu_*(\hat\xi_\lambda)=\mu(\hat\xi_{-\lambda})$. For each $z\in\C$ we can model $X$ on
the smooth completion $X_z$ of the curve
\[
X_z\setminus\{P_0,P_\infty\}\cong\{(\lambda,[v])\in\C\setminus\{0\}\times\P^2:
\hat\xi_\lambda(z) [v] = [v],\ \forall\hat\xi_\lambda\in\hat\caA\}.
\]
It is also not hard to show that all these curves are isomorphic (cf.\ \cite{FPPS})
and the completion is obtained by adding the points
$P_0=(0,[0,0,1])$ and $P_\infty=(\infty,[1,0,0])$. For each $z$ 
this embeds $X$ in $\P^1\times\P^2$ and $\caE_z$
is the pullback of the tautological bundle over $\P^2$; i.e., the bundle of common eigenlines for
all $\hat\xi_\lambda(z)$. A necessary condition for $(X,\lambda,\caL_z)$ to be the spectral
data for a minimal torus is that $\caL_z$ satisfies the reality condition
\begin{equation}
\label{eq:real}
\overline{\rho_*\caL_z}\otimes \caL_z\cong \caO_X(R)
\end{equation}
where $R$ is the divisor of ramification divisor of $\lambda$. [In fact slightly more must
be true: this linear equivalence must be achieved by a function on $X$ which is positive
over the unit circle. This is always the case if, for example,
$\caL_z\cong\caO_X(R_+)\otimes L$ where $R_+$ is the divisor of those ramification points
lying over $|\lambda|<1$ and $L$ lies in the identity component of $\{L\in
Jac(X):\overline{\rho_*L}\cong L^{-1}\}$: see \cite[\S 2.2]{McI96}.] 
\begin{lem}
For a minimal Lagrangian torus in $\CP^2$ the data $(X,\lambda,\caL_z)$ has the extra property that
\begin{equation}
\label{eq:mu}
\mu_*\caL_z\cong \caL_z^{-1}\otimes\caO_X(R).
\end{equation}
The family of line bundles $L_z=\caL_z\otimes\caL_0^{-1}$ provides a real homomorphism $L:\R^2\to
P_\R(X,\mu)$, where 
\[
P_\R(X,\mu)= \{L\in Prym(X,\mu):\overline{\rho_*L}\cong L^{-1}\}^0.
\]
This homomorphism is uniquely determined by the property that
\begin{equation}
\label{eq:Abel}
dL_0(\frac{\partial}{\partial z}) = d\caA_{P_\infty}(\frac{\partial}{\partial\zeta^{-1}})
\end{equation}
where $\caA_{P_\infty}:X\to Jac(X)$ is the Abel map with base point $P_\infty$.
\end{lem}
\begin{proof}
We will prove the equivalent statement $\mu_*\caE_z\cong \Hom(\caE_z,\caO_X(-R))$. Let $U\subset X$
be any proper open subset with $\mu(U)=U$ and $\lambda^{-1}(\lambda(U))=U$, so that
$U\to\lambda(U)$ is a 3-sheeted cover. Since $\caE\subset X\times\C^3$ any trivializing section
$v$ for $\caE_z$ over $U$ determines a matrix function $V_\lambda$ by
\[
V_{\lambda_0} = (v_{P_1},v_{P_2},v_{P_3}),\quad \lambda(P_j)=\lambda_0.
\]
Note that $\det(V_\lambda)=0$ iff $\lambda$ is a branch point. Clearly
$V_\lambda^{-1}\hat\xi_\lambda V_\lambda$ is diagonal for every $\hat\xi_\lambda\in\hat\caA$: we
will write
\[
V_\lambda^{-1}\hat\xi_\lambda V_\lambda = D(\hat\xi_\lambda).
\]
Now define $W_\lambda=\det(V_\lambda)TV_\lambda^{-1t}$, 
where $T$ is the matrix from (\ref{eq:T}).
Notice that $W_\lambda$ is holomorphic on $\lambda(U)$: it is a constant matrix times the
classical adjoint matrix for $V_\lambda$. Now we
observe that
\[
\begin{array}{rcl}
W_\lambda^{-1}\mu_*(\hat\xi_{-\lambda}) W_\lambda & = &-V_\lambda^{t}T^2\hat\xi_\lambda 
T^2V_\lambda^{-1t}\\
& = &-(V_\lambda^{-1}\hat\xi_\lambda V_\lambda)^t\\
& = & -D(\hat\xi_\lambda).
\end{array}
\]
Since $\mu:X\to X$; $(\lambda,[v])\mapsto (-\lambda,\mu[v])$ if we define
$\mu_*(\hat\xi_{-\lambda})\mu[v] = \mu[v]$ we deduce that $W_\lambda$ represents a section $w$
of $\mu_*\caE_z$ over $U$. Further, none of the columns of $W_\lambda$ vanish on $U$ 
so $w$ is a trivialising section.
On the other hand $W_\lambda^tTV_\lambda = \det(V_\lambda)$ implies that $w$ is a section of
$\Hom(T\caE_z,\caO_X(-R))$ over $U$, and $T\caE_z\cong \caE_z$ since $T$ is constant. 
Patching this argument together globally over $X$ gives the isomorphism we desire.

The linearity of $L$ and the tangent equation (\ref{eq:Abel}) are proven in \cite[Lemma 24
and p 846]{McI95}. It
must satisfy $\mu_*L\cong L^{-1}$ to respect (\ref{eq:mu}) and therefore $L$ takes values in the
Prym variety $Prym(X,\mu)=\{L\in Jac(X):\mu_*L\cong L^{-1}\}^0$.
\end{proof}
It is not hard to see that the double periodicity of the map $\varphi$ obliges $L$ to be doubly
periodic. However, it was shown in \cite[Thm 5]{McI98} that $L$ is doubly periodic 
merely when the induced metric of $\varphi$ (equally, the solution of 
(\ref{eq:tzitz})) is doubly periodic. 
As  explained in \cite[Prop 2]{McI96}, the double periodicity of the map $\varphi$ 
is equivalent to the following slightly stronger 
condition. Let $X^\prime$ be the singular curve obtained 
from $X$ by identifying the three points $O_1,O_2,O_3$ lying over $\lambda=1$. 
Let $P(X^\prime,\mu)$ be the pullback to $Prym(X,\mu)$ of the group extension 
$Jac(X^\prime)\to Jac(X)$ (whose kernel is isomorphic to $\C^*\times\C^*$)
and let 
\[
P_\R(X^\prime,\mu)=
P(X^\prime,\mu)\cap J_\R(X^\prime)^0;\quad J_\R(X^\prime)= 
\{L\in Jac(X^\prime):\overline{\rho_*L}\cong L^{-1} \}^0.
\]
Then $\varphi$ is doubly periodic iff the real homomorphism 
\begin{equation}
\label{eq:Lprime}
L^\prime:\R^2\to P_\R(X^\prime,\mu)
\end{equation}
is doubly periodic, where this is uniquely specified by the condition that
\[
dL^\prime_0(\frac{\partial}{\partial z}) = 
d\caA^\prime_{P_\infty}(\frac{\partial}{\partial\zeta^{-1}}),
\]
with $\caA^\prime_{P_\infty}:X^\prime\to Jac(X^\prime)$ the Abel map for $X^\prime$ with base
point $P_\infty$.

To understand these double periodicity conditions better, let $\G(\Omega')$ denote the space
of regular differentials on $X'$.
It is not hard to show that $P_\R(X',\mu)\simeq V/\Lambda$ where
\[
V=\{\alpha\in \G(\Omega')^*:\rho^*\alpha = -\bar\alpha,\mu^*\alpha=-\alpha\},
\]
and $\Lambda$ is the intersection of $V$ with the 
image of the lattice $H_1(X'\setminus\{O\},\Z)$ in
$\G(\Omega')^*$. Here $O$
denotes the singular point of $X'$ over $\lambda=1$. We can deduce from \cite[p 432]{McI01} 
that $V/\Lambda$ is a real compact $(p+2)$-torus whenever $X$ is smooth. 
It is easily shown that the image of $dL'$
at the origin is the real two-dimensional subspace $W\subset\G(\Omega')^*$ 
whose elements annihilate all those regular
differentials on $X'$ which vanish at $P_0$ and $P_\infty$.
Since $\Lambda$ is a full lattice, by choosing its generators we may identify $V$ with
$\R^{p+2}$. Let us use $(W,\Lambda)$ to denote the real 2-plane in $\R^{p+2}$ we obtain from
this. Then $L'$ is doubly periodic precisely when $(W,\Lambda)$ is a rational 2-plane, that
is, has a basis consisting of rational vectors.

\subsection{The moduli of minimal Lagrangian tori.}

We have learnt that every minimal Lagrangian torus gives rise to spectral 
data $(X,\lambda,\caL)$ and we know from \cite{McI95} that, conversely, any 
data satisfying the conditions described above gives rise to a minimal torus.
It remains to check that the extra symmetry of the holomorphic involution $\mu$ on $X$
forces this torus to be Lagrangian as well. This step will complete the proof of the
following theorem.
\begin{theor}
There is a bijective correspondence between:\\
(i) congruence classes of minimal Lagrangian tori $\varphi:T^2\to\CP^2$ (with base point on $T^2$), and;\\
(ii) equivalence classes of data $(X,\lambda,\caL,\mu)$ where, in addition to the reality
conditions described earlier, $\caL$ satisfies (\ref{eq:mu}) and the real homomorphism
$L^\prime$ in (\ref{eq:Lprime}) is doubly periodic.
\end{theor}
\begin{proof}
This follows from \cite{McI95,McI96,McI98} once we have established that (\ref{eq:mu}) forces
$\varphi$ to be Lagrangian. First we must recall how to construct $\varphi$ from this
spectral data. From $(X,\lambda)$ we construct the real homomorphism $L_z:\C\to
P_\R(X,\mu)$ satisfying the tangency condition (\ref{eq:Abel}). For each
$z\in\C$ we fix the line
\[
\ell_0(z) = \G(X,\caL_z(-P_0-P_\infty))\in\P\G(X,\caL_z)\cong\CP^2
\]
where $\caL_z=\caL\otimes L_z$. The reality condition (\ref{eq:real}) ensures
that $\lambda_*\caL_z$ is trivial and therefore $h^0(\caL_z)=3$. To identify
all these projective spaces to get a minimal map we need the (essentially
unique) trivialising section $\theta_z$ of $L_z$ over
$X\setminus\{P_0,P_\infty\}$ which has the properties:
(i) $\theta_z\exp(-z\zeta)$ extends holomorphically to $P_\infty$;
(ii) $\theta_z\exp(-\bar z\zeta^{-1})$ extends holomorphically to $P_0$;
(iii) $\overline{\rho_*\theta_z}=\theta_z^{-1}$,
$\mu_*\theta_z=\theta_z^{-1}$. Since $\lambda_*\caL_z$ is trivial any global
section of $\caL_z$ is determined by its restriction to the three points
$O_1,O_2,O_3$ over $\lambda=1$. So we identify $\G(X,\caL_z)$ with 
$\G(X,\caL_0)$ by mapping any global section $s$ of $\caL_z$ to the unique
global section of $\caL_0$ which has the same values as $s\theta_z^{-1}$ at
$O_1,O_2,O_3$. According to \cite{McI95,McI96} this constructs a minimal map
$\varphi:\C\to\CP^2$ which is doubly periodic precisely when the map (\ref{eq:Lprime})
is.

To show that $\varphi$ is Lagrangian it suffices to construct a loop of
frames $F_\lambda$ for which $F_1$ frames $\varphi$ and $\alpha_\lambda$
possesses the symmetry $\alpha_{-\lambda}=-T\alpha_\lambda^tT$. This frame is
constructed using the commutative diagram:
\[
\begin{array}{ccc}
\G(X\setminus\{P_0,P_\infty\},\caL_z) & \stackrel{\otimes\theta_z^{-1}}{\rightarrow}
& \G(X\setminus\{P_0,P_\infty\},\caL_0)\\
f_z\downarrow & & \downarrow f_0\\
\C(\lambda^{-1},\lambda)\otimes \C^3 & \stackrel{F_\lambda^{-1}}{\rightarrow} & 
\C(\lambda^{-1},\lambda)\otimes \C^3. 
\end{array}
\]
The maps $f_z$ are determined by fixing a unitary basis $e_j^z$
for $\G(X,\caL_z)$ (the Hermitian metric on this space 
is described in \cite[\S 3]{McI95}) for which 
\[
e_0^z\in\G(X,\caL_z(-P_0-P_\infty)),\ e_1^z\in\G(X,\caL_z(-2P_\infty)),\ 
e_2^z\in\G(X,\caL_z(-2P_0)).
\]
The scaling freedom available in the choice of this basis means $F_\lambda$ is
only determined up to right multiplication by a diagonal matrix independent of
$\lambda$. The reality conditions ensure that $
F_{\bar\lambda^{-1}}=\bar (F_\lambda^{-1})^t$ over $\vert\lambda\vert=1$ (and therefore
everywhere). Now we can use the fact that $\overline{\rho\mu_*\caL}\cong \caL$
and $\rho\mu(P_0)=P_\infty$ to insist that we also have
\[
\overline{\rho\mu_* e_0}=e_0,\quad \overline{\rho\mu_* e_1}=e_2,
\]
and therefore $F_{-\bar\lambda}^{-1}=T\bar F_\lambda T$. Putting the two
symmetries of $F_\lambda$ together shows that $F_{-\lambda}=T(F_\lambda^{-1})^tT$
whence $\alpha_{-\lambda} = -T\alpha_\lambda^tT$ as required.
\end{proof}
It is now possible to count the dimension of the moduli space of minimal Lagrangian tori in
$\CP^2$ (cf.\ \cite[\S 5]{McI98}, where this is done for minimal tori). At first we will assume that
$P_0+P_\infty$ is not a hyperelliptic divisor on $X$ (recall from \cite[\S 5]{McI98} that this case
must be treated differently, for the map $\varphi$ is $S^1$-equivariant then). 

First, the data 
$(X,\lambda,\mu)$ is completely determined by the ramification divisor $R$. This has degree
$2g+4$ but it must include the divisor $2P_0+2P_\infty$. 
Moreover, it must be preserved by the involutions $\rho$ and $\mu$,
neither of which fix any of the points in $R-2P_0-2P_\infty$, so there are $g=2p$ 
real parameters
available for $(X,\lambda,\mu)$. The line bundle $\caL$ lives in the real $p$-dimensional
space $\{\caL\otimes L:L\in P_\R(X,\mu)\}$.

Following the discussion at the end of the previous subsection, we can reduce the double 
periodicity conditions to the study of the map
\[
(X,\lambda,\mu)\mapsto (W,\Lambda)\in Gr_\R(2,\R^{p+2}).
\]
The domain and codomain both have real dimension $2p$. Given local surjectivity of this map
(the proof of which has not yet been attempted)
we will be able to find spectral data for which the plane $(W,\Lambda)$ is rational and
therefore obtain minimal Lagrangian tori. But where the map is locally surjective it must
also be injective. However, for $p>2$ this does not mean each torus will be isolated, 
since each
torus will live in a real $(p-2)$-dimensional family corresponding to the variation of $\caL$
(where we have factored out the action of translations along the universal cover of $T^2$;
that is, discarded the base point on $T^2$). Therefore we expect the moduli space to be
generically of dimension $p-2$.

In the case where $X$ is hyperelliptic with hyperelliptic divisor $P_0+P_\infty$ it follows
from \cite[\S 5]{McI98} that the map $\varphi$ is $S^1$-equivariant. Such minimal Lagrangian tori
have been thoroughly discussed \cite{CasU,Has,Joy} and in particular the spectral 
data approach
is unnecessary, since the equation (\ref{eq:tzitz}) reduces to an o.d.e.\ which can be solved
by elliptic functions. The geometrical explanation for the appearence of elliptic functions
is that, from \cite[\S 5]{McI98}, $X$
is generically genus 2 hence $Prym(X,\mu)$ is an elliptic curve. 

\subsection{The case of genus four spectral curves.}

In general the construction of examples is hampered by the difficulty in 
satisfying the periodicity conditions. Once solved, explicit formulae can 
be found using Riemann
$\theta$-functions (see, for example, \cite{Sha,McI01,MaMa}).
When $X$ has genus four the real group $P_\R(X,\mu)$ is two dimensional and therefore the
corresponding solutions of the Tzitz\'{e}ica equation are doubly periodic. 
Indeed, since $\mu$ has two fixed points $Prym(X,\mu)$ is principally polarized and therefore the
Jacobi variety of a genus two curve.  For this reason it
is interesting to describe the curves of genus four which can appear.
\begin{prop}
If $(X,\lambda)$ is the spectral data for a minimal Lagrangian torus in $\CP^2$ and $X$ has
genus $g=4$ then $(X,\lambda)$ arise from either:\\
(i) the smooth compactification of a curve with equation
\begin{equation}
\label{eq:curve1}
\lambda^2 - 2b(x) + x^3\lambda^{-2},\quad 
b(x)=b_0 + b_1x + \bar b_1 x^2 + \bar b_0x^3,\ b_j\in\C,
\end{equation}
(ii) the smooth compactification of a curve with equation
\begin{equation}
\label{eq:curve2}
k\lambda^2 - 2b(x) + k^{-1}\lambda^{-2},\quad \deg_x b(x)=3,\ b(x)\in\R[x],\ |k|=1.
\end{equation}
\end{prop}
\begin{proof}
From our discussion above we see that the covering $\delta:X\to Y$ 
uniquely determines $X$ as the covering curve for 
$\lambda=\sqrt{y}$ where
$y:Y\to\hat\C$ is the rational function of degree three with a single zero at
$Q_0=\delta(P_0)$ and a single pole at $Q_\infty=\delta(P_\infty)$. 
All we require of $(Y,y)$ is that it have genus 2 and 
possess a real involution $\rho$ for which $\overline{\rho^*y}=y^{-1}$ and which fixes
every point over $|y|=1$, on which there are no branch points. Now we recall that
every genus 2 curve has a unique hyperelliptic involution: let $\iota\in Aut(X)$
denote this involution. Since $\rho\iota\rho$ is also a hyperelliptic involution we
deduce that $\rho\iota=\iota\rho$. Now we separate the argument into two cases: (i)
$\iota(Q_0)\neq Q_\infty$, (ii) $\iota(Q_0)=Q_\infty$.

In case (i) we can choose degree 2 function $x$ with divisor
$(x)=Q_0+\iota(Q_0)-Q_\infty-\iota(Q_\infty)$ and this therefore has $x\overline{\rho^*x}$
constant. We can rescale by a positive constant to make $x\overline{\rho^*x}=1$. By
comparing divisors it follows that $y\iota^*y x^{-3}$ is also constant and by a
unimodular scaling of $x$ we may assume this is positive, $y\iota^*y x^{-3}=a>0$. 
But
\[
a^{-3}x^{-3}=\overline{\rho^*(y\iota^*y)}=a^3x^{-3},
\]
so $a=1$. Thus $y\iota^*y =x^3$. Again by considering divisors we deduce that
$y+\iota^*y=2b(x)$ where $b(x)$ is a polynomial of degree 3. Therefore 
\[
y - 2b(x) +x^3y^{-1}=0.
\]
From this we compute $\overline{b(\bar x)}=x^3b(x^{-1})$ so $b(x)$ has the form 
given in (\ref{eq:curve1}).

In case (ii) we choose $x$ to have divisor of poles $Q_0+Q_\infty$ and so that
$\overline{\rho^*x}=x$. In this case $y\iota^*y$ is a unimodular constant $k^2$. Set $\tilde y
= y/k$ then $\tilde y\iota^*\tilde y=1$. As above $\tilde y +\iota^*\tilde y = 2b(x)$
for some polynomial of degree 3. Therefore
\[
ky - 2b(x) + k^{-1}y^{-1}=0,\ |k|=1,
\]
and a computation shows that $b(x)\in\R[x]$.
\end{proof}
This description does not include the conditions on the ramification divisor
$R$, which are somewhat complicated (however, see the Remark below).
In standard hyperelliptic form the corresponding curves $Y$ are as follows. For
case (i) set $w=y-b(x)$ then
\[
w^2 = b(x)^2-x^3.
\]
For case (ii) set $w=ky-b(x)$, then
\[
w^2 = b(x)^2-1.
\]
It is possible to identify $Prym(X,\mu)$ with the Jacobian of a genus two
curve lying in the symmetric product $S^2X$, using Donagi's bigonal construction
\cite{Don}. However, let us finish by simply pointing out that for curves of
the type (\ref{eq:curve2}) the Prym variety is isogenous to a product of
elliptic curves. This echoes the construction of minimal Lagrangian tori by
Joyce in \cite[\S 6.2]{Joy}, where he finds formulae which involve products of
two different elliptic functions. The spectral curves he produces for these examples are all
of type (\ref{eq:curve2}).
\begin{prop}
Let $E_1,E_2$ be the elliptic curves with respective equations 
\[
(i)\ w^2=b(x)+1;\quad (ii)\ w^2=b(x)-1,\quad \deg_x b(x)=3,\ b(x)\in\R[x].
\]
Then for $X$ of type (\ref{eq:curve2}), $Prym(X,\mu)$ is isogenous to
$E_1\times E_2$.
\end{prop}
\begin{proof}
Any curve $X$ with equation (\ref{eq:curve2}) admits two involutions
$\iota_1,\iota_2$ given by
\[
\iota_1(x,\lambda)=(x,\lambda^{-1});\quad \iota_2(x,-\lambda^{-1}).
\]
Notice that $\iota_1\iota_2=\mu$. The respective rings of invariants are
\[
\C[x,\lambda]^{\iota_1}=\C[x,\lambda+\lambda^{-1}],\quad
\C[x,\lambda]^{\iota_2}=\C[x,\lambda-\lambda^{-1}].
\]
Taking, respectively, $w=(\lambda+\lambda^{-1})/\sqrt{2}$ 
and $w=(\lambda-\lambda^{-1})/\sqrt{2}$ yields the equations 
for $E_1$ and $E_2$. Let $f_j:X\to E_j$ be the corresponding coverings, then
we have the isogeny 
\[
Jac(X)\approx f_1^*E_1\times P_1\approx f_2^*E_2\times P_2
\]
where $P_1,P_2$ are the respective Prymians. Since $\iota_1,\iota_2$ commute
it follows that $f_1^*E_1\cap P_2$ and $f_2^*E_2\cap P_1$ are distinct one
dimensional subgroups. Now if $L\in f_1^*E_1\cap P_2$ then
\[
\mu^*L=\iota_1^*\iota_2^*L\cong \iota_1^*L^{-1}\cong L^{-1},
\]
and similarly for $L\in f_2^*E_2\cap P_1$. It follows by considering
dimensions that
\[
Prym(X,\mu)\approx (f_1^*E_1\cap P_2)\times (f_2^*E_2\cap P_1)\approx
E_1\times E_2.
\]
\end{proof}
\textit{Remark.} There is a simple characterization of
the polynomials $b(x)\in\R[x]$ which give
spectral data of the type (\ref{eq:curve2}); i.e., those polynomials
for which the real involution satisfies the right conditions over $|\lambda|=1$. 
\begin{lem}
A pair $(X,\lambda)$ with equation (\ref{eq:curve2}) satisfies all the
properties for being spectral data if and only if for some scaling of $x$
\[
b(x) = x^3 - \frac{3}{2} (u+v)x^2+3uvx+w,\ u,v,w\in\R,
\]
where $v-u>4^{1/3}$ and
\[
1+\frac{1}{2}u^2(u-3v)<w<\frac{1}{2}v^2(v-3u)-1.
\]
\end{lem}
\begin{proof}
By rescaling $x$ we may assume $b(x)$ is monic.
We must show that these are precisely the conditions under which
$\lambda:X\to\hat\C$ has no branch point on $|\lambda|=1$ and
$\rho(x,\lambda)=(\bar x,\bar\lambda^{-1})$ fixes every point over
$|\lambda|=1$. The locus of points on $X$ over the unit circle is
characterized by the equation $b(x)-\cos(\theta)=0$, $\theta\in[0,2\pi]$: we
see this by setting $\lambda^2k=e^{i\theta}$ in (\ref{eq:curve2}). We require
that this equation has distinct real roots for all values of $\theta$. If we
write
\[
b^\prime(x)=3(x-u)(x-v),\ u<v,
\]
then since $|\cos(\theta)|<1$ we require $b(u)>1$ and $b(v)<-1$: these two
give the last inequality in the lemma. In particular
$2<b(u)-b(v)=(v-u)^3/2$, which gives $v-u>4^{1/3}$. 
\end{proof}

\subsection{Totally real superconformal maps into $\CP^{2k}$.} 

When $n$ is even those superconformal maps $\varphi:M\to\CP^n$ which are
``totally real'' (in the sense that $\varphi^*\omega = 0$, see the remark below)
can also be characterized by the existence of a certain primitive
harmonic lift. Essentially by definition \cite{BolPW},
superconformal maps have an orthogonally periodic harmonic
sequence and therefore admit a primitive lift $\psi:M\to SU_{n+1}/H$, where $H$ is the
maximal torus in $SU_{n+1}$: the target is the full flag manifold over $\CP^n$. 
Then $\varphi$ is totally real
if and only if it admits a further primitive lift
$\tilde\psi:\tilde M\to SU_{n+1}/H_0$ where $H_0\subset H$ is the $n/2$-dimensional subgroup
fixed by the outer involution $\mu$ of $\fa_n$ induced by the $\Z_2$ symmetry of its Dynkin diagram.
To see this, write $k=n/2$ and order the harmonic sequence $\ell_0,\ell_1,\ldots,\ell_n$ so
that $\varphi$ corresponds to $\ell_k$. Following \cite{BolW} we choose local sections $f_j$
of $\ell_j$ so that
\[
\frac{\partial f_j}{\partial z} = f_{j+1}+\frac{\partial}{\partial z}\log |f_j|^2 f_j,\
j=0,\ldots n-1.
\]
These imply the Toda equations
\[
\frac{\partial^2}{\partial z\partial\bar z}\log(|f_j|^2) = \frac{|f_{j+1}|^2}{|f_j|^2} -
\frac{|f_j|^2}{|f_{j-1}|^2},\ j=0,\ldots,n-1.
\]
Since $\varphi$ is totally real
we can arrange for $|f_k|=1$. An inductive argument using the
Toda equations then shows that $|f_{k+j}|=|f_{k-j}|^{-1}$ for $j=0,\ldots,k$. 
Now let $F:\tilde M\to SU_{n+1}$ be the Toda frame, whose $j$-th column is $f_j/|f_j|$, and
let $\alpha = F^{-1}dF$ be its Maurer-Cartan form.
Then a simple computation shows that $\alpha_\fh$ lies in $\fg^\nu_0\cap\fg^\mu_0$ and
$\alpha_\fp^{1,0}$ lies in $\fg^\nu_1\cap\fg^\mu_1$, where $\nu$ is again the Coxeter-Killing
automorphism. Therefore, with an appropriate initial condition for $F$, all the maps in this
harmonic sequence have a (local) primitive lift into $SU_{n+1}/H_0$, since this is the
$2n+2$-symmetric space for the product automorphism $\nu\mu$. 

In particular, for totally real superconformal tori we can argue as above that the lift 
factors through a covering torus, with covering group a subgroup of $\Z_{n+1}$. The spectral
data for such tori will consist of $(X,\lambda,\caL,\mu)$ where $\mu$ is a holomorphic
involution which the other data respects in the same fashion as for the case $n=2$. 

\smallskip\noindent
\textit{Remark.} There is a schism in the mathematical community as to the
definition of a totally real submanifold. It seems that the original meaning,
used in the 1960's in the study of real submanifolds $M$ of a complex 
submanifold, is that $TM\cap J(TM)=0$, where $J$ is the complex structure on the
ambient space. Somewhere in the 1970's this term started to 
be used also for the condition that $TM$ is perpendicular to $J(TM)$, which is the
sense I have used above. Unfortunately, we now have the situation where both 
definitions are in use and it is not always clearly stated which definition is meant.

\end{document}